\documentclass{amsart}
\usepackage{amssymb}
\usepackage{amsmath}
\usepackage{amsthm}
\usepackage{mathrsfs}
\usepackage{colortbl}
\usepackage{graphicx}

\theoremstyle{definition}

\newtheorem{theoremz}{Theorem}[section]
\newtheorem{corollaryz}[theoremz]{Corollary}
\newtheorem{lemmaz}[theoremz]{Lemma}

\newtheorem{definitionz}[theoremz]{Definition}
\newtheorem{assumptionz}[theoremz]{Assumption}
\newtheorem{remarkz}[theoremz]{Remark}

\makeatletter

\@addtoreset{theorem}{section}
\makeatother

\makeatletter

\@addtoreset{equation}{section}
\makeatother

\begin{document}                                                 
\title[Error Analysis]{Error Analysis of Approximate Operators for a Particle Method based on Voronoi Diagram}                                 
\author[Hajime Koba]{Hajime Koba}                                
\address{Graduate School of Engineering Science, Osaka University,\\
1-3 Machikaneyamacho, Toyonaka, Osaka, 560-8531, Japan}                                  
\email{iti@sigmath.es.osaka-u.ac.jp}

\author[Kazuki Sato]{Kazuki Sato}                                
\address{Graduate School of Engineering Science, Osaka University,\\
1-3 Machikaneyamacho, Toyonaka, Osaka, 560-8531, Japan}                                  
\email{k-sato@sigmath.es.osaka-u.ac.jp}

\begin{abstract}
This paper considers several approximate operators used in a particle method based on a Voronoi diagram. We introduce and study our approximate operators on gradient and Laplace operators. We derive error estimates for these approximate operators by applying our weight functions. The key idea of deriving our error estimates is to divide the integration region into a ring-shaped area and some areas. In the Appendix, we give an example application of the main results of this paper.
\end{abstract}       
\maketitle

\section{Introduction}\label{sect1}

\begin{figure}[htbp]
\includegraphics[width=12cm]{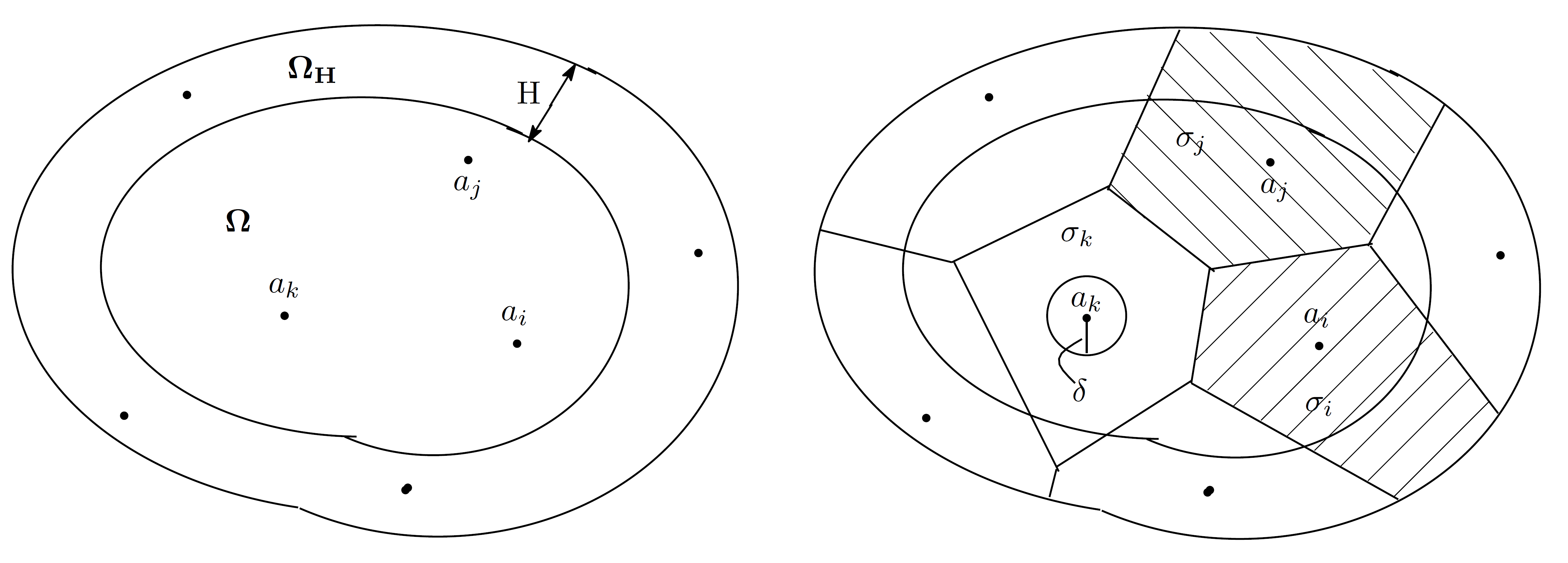}
\caption{Arrangement of Particles and Voronoi Diagram}
\label{Fig1}
\end{figure}

We are interested in the error analysis of approximate operators for a moving particle semi-implicit method (MPS). A moving particle semi-implicit method is a numerical method developed by Koshizuka-Oka \cite{KO96}. In \cite{KO96}, they introduced approximate operators of gradient and Laplace operators based on a Voronoi diagram. Ishijima-Kimura \cite{IK10} considered the error on the approximate operators under assumptions on their weight function. In \cite{IK10}, they applied the spherical symmetry of their weight function to derive their error estimates. Imoto-Tagami \cite{IT16, IT17} modified the approximate operators introduced in \cite{KO96}, and studied their error estimates for their approximate operators. In \cite{IT16,IT17}, they derived their convergence rates of their error estimates with respect to influence radius by using assumptions on their weight function.

This paper has three purposes. The first one is to generalize the approximate operators introduced in \cite{KO96,IT16,IT17}. More precisely, this paper standardizes their approximate operators (see Definition \ref{def12} for details). The second is to derive error estimates for our approximate operators. We use some properties of our weight functions (Assumption \ref{ass11}) to study the error on our operators. The key idea of deriving our error estimates is to divide the integration region into a ring-shaped area and some areas. The third is to give an application of our main results (see the Appendix for details).

Let us first introduce notations. Let $x = { }^t ( x_1 , x_2), y = { }^t ( y_1 , y_2) , z = { }^t (z_1, z_2) \in \mathbb{R}^2$ be the spatial variables, and let $\alpha = ( \alpha_1 , \alpha_2) \in \mathbb{N}_0^2$ be a multi-index, where $\mathbb{N}_0 := \mathbb{N} \cup \{ 0 \}$. For each multi-index $\alpha = ( \alpha_1 , \alpha_2 )$,
$| \alpha | := \alpha_1 + \alpha_2$, $D^\alpha := \partial_1^{\alpha_1} \partial_2^{\alpha_2}$, and $\alpha ! := \alpha_1 ! \alpha_2 !$, where $\partial_j := \partial/{\partial x_j}$, $\partial_j^0 : = 1$, and $0 ! :=1$. For each $x \in \mathbb{R}^2$, $r , \ell >0$, the symbols $B_r$, $B_r (x)$, $B_{r , \ell}$, and $B_{r , \ell} (x)$ are defined by $B_r = \{ y \in \mathbb{R}^2; { \ } | y | < r \}$, $B_r (x) = \{ y \in \mathbb{R}^2; { \ } | x - y | < r \}$,
\begin{equation*}
B_{r , \ell} = \{ y \in \mathbb{R}^2;{ \ } \ell < | y | < r \}, \text{ and }B_{r , \ell} (x) = \{ y \in \mathbb{R}^2;{ \ } \ell < | y - x | < r \}. 
\end{equation*}
Note that $B_{r, \ell} = B_r \setminus \overline{B}_\ell$. Let $U \subset \mathbb{R}^n$ be a domain, and $\overline{U}$ be the closure of $U$. For each smooth function $f = f(x)$ on $\overline{U}$, we define
\begin{align*}
& \| f \|_{C^m} = \| f \|_{C^m ( \overline{U})} = \max_{ j \in \{ 1,2, \cdots , m \}} | f |_{C^j ( \overline{U})},{ \ \ }\| f \|_{L^1} = \| f \|_{L^1 (U)} = \int_U | f (x)| { \ } d x,\\
& |  f |_{C^0} = | f |_{C^0 ( \overline{U})} = \max_{x \in \overline{U}} | f (x)|,{ \ \ }| f |_{C^j} = | f |_{C^j ( \overline{U})} = \max_{ | \alpha | = j } | D^\alpha f |_{C^0 ( \overline{U})},\\
& \| { \ } | \cdot |  f  (  \cdot ) \|_{L^1 (U)} = \int_U |x | |f (x)| { \ } dx,{ \ \ }\| { \ } | \cdot |^2 f (  \cdot ) \|_{L^1 (U)} = \int_U |x|^2 | f (x)|^2  { \ } dx.
\end{align*}
Morevoer, for each $g \in L^1 ( B_{r , \ell} (a)) $ for some $a \in \mathbb{R}^2$ and $s \in \mathbb{R}$, we define
\begin{equation*}
\| { \ } | a - \cdot |^s  g  ( a -  \cdot ) \|_{L^1 (B_{r, \ell} (a))} = \int_{B_{r, \ell} (a)} | a - y |^s |g ( a - y)| { \ } d y.
\end{equation*}

The following notations are of particular importance in this study. Let $\Omega \subset \mathbb{R}^2$ be a bounded domain, and $H$ be a positive constant. Throughout this paper, we fix $\Omega $ and $H$. Define
\begin{equation*}
\Omega_H = \{ x \in \mathbb{R}^2; { \ } | x - y | < H, { \ }y \in \Omega  \}.
\end{equation*}
See Fig. \ref{Fig1}. Let $N \in \mathbb{N}$ and $a_1 , a_2 , \cdots , a_N \in \Omega_H$ such that $a_i \neq a_j$ if $i \neq j$. Write
\begin{equation*}
\Omega_{N, H} = \{ a_i \}_{ i =1}^N .
\end{equation*}
For each $i \in \{ 1, 2 , \cdots , N \}$, we define
\begin{equation*}
\sigma_i := \{ x \in \Omega_H; { \ }| x - a_i | < | x - a_j | \text{ for each } a_j \in \Omega_{N,H} { \ }(j \neq i ) \} .
\end{equation*}
See Fig. \ref{Fig1}. In general, we call $\sigma_i$ a \emph{Voronoi region}, and $ \{ \sigma_i \}_{i = 1}^N$ a \emph{Voronoi diagram} (see \cite{Vor1908}). Since
\begin{equation*}
\overline{\Omega}_H = \bigcup_{ i = 1}^N \overline{\sigma}_i \text{ and } \sigma_i \cap \sigma_j = \emptyset { \ }(i \neq j),
\end{equation*}
we call $\sigma_i$ a \emph{Voronoi cell} and $\{ \sigma_i \}_{i = 1}^N$ a \emph{Voronoi decomposition}. Write
\begin{equation*}
r_\sigma = \max_{i \in \{ 1 , \cdots, N \}} \max_{y \in \overline{\sigma}_i}\{ | y - a_i | \}.
\end{equation*}
By definition, we see that $\sigma_i \subset \overline{B}_{r_\sigma} ( a_i )$. Throughout this paper, we assume that
\begin{equation*}
r_\sigma < H .
\end{equation*}
Let $h >0$ such that $r_\sigma < h < H$. Let $k \in \{1 , \cdots, N \}$ be such that
\begin{equation*}
a_k \in \sigma_k \cap \Omega.
\end{equation*}
Let $\delta > 0$ such that $B_\delta (a_k) \subset \sigma_k \cap \Omega$ (see Fig. \ref{Fig1}). Throughout this paper, we fix $r_\sigma$, $h$, $k$, and $\delta$. For each $x \in B_\delta (a_k)$,
\begin{align*}
\overline{\mathcal{R}} ( x, h ) & := \{ i \in \{ 1 , \cdots, N \} ; { \ } 0 \leq | x - a_i | < h \},\\
\mathcal{R}( x, h ) & := \{ i \in \{ 1 , \cdots , N \} ; { \ } 0 \leq | x - a_i | < h \} \setminus \{ k \},
\end{align*}
and for each $i \in \{ 1 , \cdots , N \}$,
\begin{equation*}
V_i (x) : = \int_{\sigma_i \cap B_{h , \delta } (x)} 1 { \ }d x .
\end{equation*}
We assume that $B_{h + r_\sigma} (x) \subset \Omega_H$ for each $x \in B_\delta (a_k)$. We call $h$ the \emph{radius of the interaction area (influence radius)}.

Next we introduce the assumptions of our weight functions and approximate operators.
\begin{assumptionz}[Weight functions]\label{ass11}
Let $w \in L^\infty ( \mathbb{R}^2 ) \cap C ( \overline{ B}_{h , \delta} )$. We call $w$ a \emph{weight function} if the following five properties hold:\\ 
$(\mathrm{i})$ For almost all $x \in \mathbb{R}^2$,
\begin{equation*}
w (x) \geq 0.
\end{equation*}
$(\mathrm{ii})$ For almost all $x \in \mathbb{R}^2 \setminus \overline{B}_{h , \delta}$,
\begin{equation*}
 w (x) = 0.
\end{equation*}
$(\mathrm{iii})$ There is $L_w \geq 0$ such that for all $x , y \in B_{h , \delta }$
\begin{equation*}
| w (x) - w (y) | \leq L_w | x - y |.
\end{equation*}
$(\mathrm{iv})$ There is $\widehat{w} \in L^1 ( \delta , h )$ such that for almost all $x \in B_{h , \delta}$
\begin{equation*}
w (x) = \widehat{w}(| x |).
\end{equation*}
$(\mathrm{v})$ There is $C_0 >0$ such that for each $x \in B_{\delta} (a_k)$,
\begin{align*}
\sum_{i \in \mathcal{R} ( x ,h) } \int_{\sigma_i} w ( x - y) { \ } d y & \geq C_0,\\
\sum_{j \in \mathcal{R} (x ,h) } V_j ( x ) w ( x - a_j ) & \geq C_0. 
\end{align*}
\end{assumptionz}
Remark that several radial functions are our weight functions. In fact, we set
\begin{equation*}
w (x) = 
\begin{cases}
|x|^m , & x \in \overline{B}_{h, \delta},\\
0, & x \in \mathbb{R}^2 \setminus B_{h, \delta},
\end{cases}
\end{equation*}
for some $m \in \mathbb{Z}$. Since $|x | - |y| = (|x|^2 - |y|^2)/(|x| + |y|) \leq C(h, \delta ) | x - y|$ for $x ,y \in \overline{B}_{h, \delta}$, we easily check that $w$ satisfies the properties of Assumption \ref{ass11}. Applying the spherically symmetric property of our weight function, we derive our error estimates. See Sections \ref{sect3}-\ref{sect6} for details.

\begin{definitionz}[Approximate operators]\label{def12}
Let $w$ be a weight function satisfying the properties of Assumption \ref{ass11}. For each $f \in C ( \overline{\Omega}_H )$, define the approximate operators as follows:
\begin{equation*}
\widetilde{\Pi}_h f ( a_k ) = \frac{ \displaystyle{ \sum_{ i \in \mathcal{R} ( a_k , h)} V_i ( a_k ) f ( a_i ) w ( a_k - a_i )} }{  \displaystyle{ \sum_{j \in \mathcal{R} ( a_k , h )} V_j ( a_k )  w ( a_k - a_j )  }  },
\end{equation*}
\begin{equation*}
\widetilde{\nabla}_h f ( a_k ) = 2 \frac{  \displaystyle{ \sum_{i \in \mathcal{R} ( a_k , h)} V_i ( a_k ) \frac{f (a_k) - f (a_i)}{| a_k - a_i|} \frac{ a_k - a_i }{| a_k -a_i|} w ( a_k - a_i ) }  }{  \displaystyle{ \sum_{j \in \mathcal{R} ( a_k ,h)} V_j ( a_k ) w ( a_k - a_j )  }  },
\end{equation*}
\begin{equation*}
\widetilde{\Delta}_h f ( a_k ) = -4 \frac{  \displaystyle{  \sum_{i \in \mathcal{R} ( a_k , h)} V_i ( a_k ) \{ f ( a_k ) - f (a_i) \} w ( a_k - a_i ) }  }{  \displaystyle{ \sum_{j \in \mathcal{R} ( a_k , h )} V_j ( a_k ) | a_k  - a_j |^2 w ( a_k - a_j )  }  },
\end{equation*}
\begin{equation*}
\widetilde{\square}_h f ( a_k ) = -4 \frac{  \displaystyle{  \sum_{i \in \mathcal{R} ( a_k , h)} V_i ( a_k ) \frac{ f ( a_k ) - f (a_i) }{ | a_k - a_i |^2 } w ( a_k - a_i ) }  }{  \displaystyle{ \sum_{j \in \mathcal{R} ( a_k , h )} V_j ( a_k )  w ( a_k - a_j )  }  }.
\end{equation*}
\end{definitionz}
Note that $| \widetilde{\Pi}_h f ( a_k ) | \leq | f |_{C^0} $, $| \widetilde{\nabla}_h f ( a_k ) | \leq 4 \delta^{-1} | f |_{C^0} $, $| \widetilde{\Delta}_h f ( a_k ) | \leq 8 \delta^{-2} | f |_{C^0} $, and $| \widetilde{\square}_h f ( a_k ) | \leq 8 \delta^{-2} | f |_{C^0} $.

The main results of this paper are as follows.
\begin{theoremz}\label{thm13}
Let $w$ be a weight function satisfying the properties of Assumption \ref{ass11}. Then for each $f \in C^1 (\overline{\Omega}_H)$,
\begin{equation*}
| f ( a_k ) - \widetilde{\Pi}_h f ( a_k ) | \leq ( h + r_\sigma ) | f |_{C^1} + \{ 2 c_1 (a_k) + 2 c_2 (a_k) \} | f |_{C^0}.
\end{equation*}
Here
\begin{align*}
c_1 (a_k ) & := \frac{ \| w (a_k - \cdot ) \|_{L^1 (\sigma_k \setminus B_\delta (a_k) )}}{ \| w ( a_k - \cdot ) \|_{L^1 (B_h (a_k) \setminus B_\delta (a_k) )}},\\
c_2 (a_k) & := \frac{ \pi L_w r_\sigma h^2 }{ \displaystyle{ \| w ( a_k - \cdot ) \|_{L^1 ( B_h ( a_k ) \setminus \sigma_k )} } }.
\end{align*}
\end{theoremz}

\begin{theoremz}\label{thm14}
Let $w$ be a weight function satisfying the properties of Assumption \ref{ass11}. Assume that $\mathcal{R}( a_k , \lambda h ) \neq \emptyset$ for some $0 < \lambda < 1$. Then for each $f \in C^2 (\overline{\Omega}_H)$,
\begin{equation*}
| \nabla f ( a_k ) - \widetilde{\nabla}_h f (a_k ) | \leq 4 h | f |_{C^2} + \left\{ \frac{8 r_\sigma}{ \lambda h} + 4 c_1 (a_k) + 4 c_2 (a_k) + 8 c_3 ( a_k ) \right\} | f |_{C^1}.
\end{equation*}
Here $c_1 (a_k)$, $c_2 (a_k)$ are the constants defined by Theorem \ref{thm13}, and
\begin{equation*}
c_3 ( a_k) : = \frac{ \displaystyle{ \| w ( a_k - \cdot ) \|_{L^1 ( B_{ \lambda h + r_\sigma } ( a_k ) \setminus \sigma_k )} } }{\| w ( a_k - \cdot ) \|_{L^1 ( B_h ( a_k ) \setminus \sigma_k )}} .
\end{equation*}
\end{theoremz}

\begin{theoremz}\label{thm15}
Let $w$ be a weight function satisfying the properties of Assumption \ref{ass11}. Assume that $\mathcal{R}( a_k , \lambda h ) \neq \emptyset$ for some $0 < \lambda < 1$. Then for each $f \in C^3 (\overline{\Omega}_H)$,
\begin{equation*}
| \Delta f ( a_k ) - \widetilde{\Delta}_h f ( a_k ) | \leq 24 h | f |_{C^3} + \left\{ 4 \sum_{ \mathfrak{i}=4}^7 c_{\mathfrak{i}} (a_k) \right\} | f |_{C^1}.
\end{equation*}
Here
\begin{multline*}
c_4 (a_k) : = \frac{ \| | a_k - \cdot | w (a_k - \cdot ) \|_{L^1 (\sigma_k \setminus B_\delta (a_k) ) } }{ \| | a_k - \cdot |^2 w (a_k - \cdot ) \|_{L^1 (B_h (a_k) \setminus B_\delta (a_k) )} }\cdot\\ 
\bigg( 1 + r_\sigma \frac{ \| | a_k - \cdot | w (a_k - \cdot ) \|_{L^1 ( B_h (a_k) \setminus B_\delta (a_k) ) } }{ \| | a_k - \cdot |^2 w (a_k - \cdot ) \|_{L^1 (B_h (a_k) \setminus \sigma_k )} } \bigg),
\end{multline*}
\begin{align*}
c_5 (a_k) & := \frac{ r_\sigma }{ \lambda h} \frac{ \displaystyle{ \| | a_k - \cdot | w (a_k - \cdot ) \|_{L^1 ( B_h (a_k) \setminus \sigma_k )} }  }{ \displaystyle{ \| | a_k - \cdot |^2 w ( a_k - \cdot ) \|_{L^1 (B_h (a_k) \setminus \sigma_k ) } }  },\\
c_6 (a_k) & := r_\sigma \frac{ \displaystyle{ \| w (a_k - \cdot ) \|_{L^1 ( B_{\lambda h + r_\sigma} (a_k) \setminus \sigma_k )} }  }{ \displaystyle{ \| | a_k - \cdot |^2 w ( a_k - \cdot ) \|_{L^1 (B_h (a_k) \setminus \sigma_k ) } }  },\\
c_7 ( a_k ) & := \frac{2 r_\sigma h \| w (a_k - \cdot ) \|_{L^1 (B_h (a_k) \setminus \sigma_k )} }{ \| | a_k - \cdot |^2 w (a_k - \cdot ) \|_{L^1 (B_h (a_k) \setminus \sigma_k )}} \frac{ \displaystyle{\sum_{i \in \mathcal{R} (a_k , h )} V_i | a_k - a_i | w (a_k - a_i)}}{ \displaystyle{\sum_{j \in \mathcal{R} (a_k , h )} V_j | a_k - a_j |^2 w (a_k - a_j)  } } ,
\end{align*}
and
\begin{equation*}
c_8 ( a_k ) := \frac{ \pi L_w r_\sigma h^3 }{ \| | a_k - \cdot |^2 w ( a_k - \cdot ) \|_{L^1 (B_h (a_k) \setminus \sigma_k )} } \Bigg( 1 + h \frac{ \displaystyle{\sum_{i \in \mathcal{R} (a_k , h )} V_i | a_k - a_i | w (a_k - a_i)}}{ \displaystyle{\sum_{j \in \mathcal{R} (a_k , h )} V_j | a_k - a_j |^2 w (a_k - a_j)  } }  \Bigg).
\end{equation*}

\end{theoremz}

\begin{theoremz}\label{thm16}
Let $w$ be a weight function satisfying the properties of Assumption \ref{ass11}. Assume that $\mathcal{R}( a_k , \lambda h ) \neq \emptyset$ for some $0 < \lambda < 1$. Then for each $f \in C^3 (\overline{\Omega}_H)$,
\begin{equation*}
| \Delta f ( a_k ) - \widetilde{\square}_h f ( a_k ) | \leq 24 h | f |_{C^3} + \{ 8 c_9 (a_k) + 16 c_{10} (a_k) + 16 c_{11} (a_k) + 4 c_{12} (a_k) \} | f |_{C^1}.
\end{equation*}
Here
\begin{align*}
c_9 (a_k) & := \frac{ \displaystyle{ \|  w ( a_k - \cdot )/| a_k - \cdot | \|_{L^1 ( \sigma_k \setminus B_\delta (a_k) )} }  }{\| w ( a_k - \cdot ) \|_{L^1 ( B_h ( a_k ) \setminus B_\delta (a_k) )}} ,\\
c_{10} (a_k) & := \frac{r_\sigma}{\lambda h} \frac{ \displaystyle{ \left\| w ( a_k - \cdot )/| a_k - \cdot | \right\|_{L^1 ( B_h (a_k) \setminus \sigma_k )} }  }{ \displaystyle{ \| w ( a_k - \cdot ) \|_{L^1 (B_h ( a_k ) \setminus \sigma_k )} }  },\\
c_{11} (a_k) & := \frac{ \displaystyle{ \left\| w ( a_k - \cdot )/| a_k - \cdot | \right\|_{L^1 ( B_{\lambda h + r_\sigma } (a_k) \setminus \sigma_k )} }  }{ \displaystyle{ \| w ( a_k - \cdot ) \|_{L^1 (B_h ( a_k ) \setminus \sigma_k )} }  },
\end{align*}
and
\begin{equation*}
c_{12} (a_k) : = \frac{ \pi L_w }{ \| w (  a_k  - \cdot ) \|_{L^1 ( B_h ( a_k ) \setminus \sigma_k )} } \Bigg( \frac{ 2 r_\sigma h}{ \lambda} + ( \lambda h + r_\sigma )^2 + \frac{ \displaystyle{\sum_{i \in \mathcal{R} (a_k , \lambda h )} V_i \frac{w (a_k - a_i)}{|a_k - a_i|} } }{ \displaystyle{\sum_{j \in \mathcal{R} (a_k , h )} V_j w (a_k - a_j)  } }  \Bigg).
\end{equation*}
\end{theoremz}

\begin{remarkz}\label{rem17}
We explain the application of our results to numerical analysis and simulation. In (MPS), $\Omega_{N,H} = \{ a_i \}_{i=1}^N$ means the distribution of the particles in the domain $\Omega_H$, $a_i$ the position of a particle, and $N$ the number of the particles in $\Omega_H$. Using our approximate operators, we can study an approximation of a fluid system, however, we need advanced techniques.  From Theorems \ref{thm13}-\ref{thm16}, we see that our approximate operators become better approximations when the number of particles is sufficiently large. 
\end{remarkz}

In Section \ref{sect2}, we study some properties of our weight functions. We prove Theorem \ref{thm13} in Section \ref{sect3}, Theorem \ref{thm14} in Section \ref{sect4}, Theorem \ref{thm15} in Section \ref{sect5}, and Theorem \ref{thm16} in Section \ref{sect6}. In the Appendix, we give an application of the main results of this paper.

\section{Preliminaries}\label{sect2}

In this section, we recall the Taylor theorem and study some fundamental properties of our weight functions.
\begin{lemmaz}[Taylor's theorem]\label{lem21}
Let $m \in \mathbb{N}$ and $f \in C^{ m + 1} ( \overline{\Omega}_H)$. Then for each $x , y \in \Omega_H$,
\begin{equation*}
f (y) = f (x) + \sum_{ 1 \leq | \alpha | \leq m} \frac{D^\alpha f (x)}{\alpha !} ( y - x )^\alpha + R_{m + 1}[f] (x , y ). 
\end{equation*}
Here
\begin{equation*}
R_{ m + 1 }[f] (x,y) := \sum_{| \alpha | = m + 1} (y - x )^\alpha \frac{m + 1}{\alpha !} \int_0^1 (1 - t)^m D^\alpha f ( t y + (1 - t ) x ) { \ } d t.
\end{equation*}
Moreover, for each $x , y \in \Omega_H$,
\begin{equation}
| R_{ m + 1 }[f] (x,y) | \leq 2 ( m + 1 )| x - y |^{ m + 1 } | f |_{C^{m + 1 }} . \label{eq21}
\end{equation}
\end{lemmaz}

\begin{proof}(Lemma \ref{lem21})
We only derive \eqref{eq21}. To this end, we show that for each $m \in \mathbb{N}$
\begin{equation}
\sum_{| \alpha | = m} \frac{1}{\alpha !} \leq 2. \label{eq22}
\end{equation}
A direct calculation gives
\begin{align*}
\sum_{| \alpha | = 1} \frac{1}{\alpha !} & = \frac{1}{1 ! 0 !} + \frac{1}{ 0! 1 !} = 1 + 1 =2,\\
\sum_{| \alpha | = 2} \frac{1}{\alpha !} & = \frac{1}{2! 0 !} + \frac{1}{ 1! 1 !} + \frac{1}{0 ! 2 !} = \frac{1}{2} + 1  + \frac{1}{2}=2.
\end{align*}
We now assume that $m \geq 3$. Let $\alpha = (\alpha_1 , \alpha_2) \in \mathbb{N}_0^2$ such that $\alpha_1 + \alpha_2 = m$. It is easy to check that
\begin{equation*}
\frac{1}{\alpha !} = \frac{1}{\alpha_1 ! \alpha_2 !} \leq \left( \frac{1}{2} \right)^{\alpha_1 - 1}\left( \frac{1}{2} \right)^{\alpha_2 - 1} =  \left( \frac{1}{2} \right)^{m - 2}{ \ \ } (\text{whenever } \alpha_1 \alpha_2 \neq 0)
\end{equation*}
and that
\begin{equation*}
\frac{1}{\alpha_1 !} \leq \left( \frac{1}{2} \right)^{m - 1} \leq  \left( \frac{1}{2} \right)^{m - 2} \text{ and  } \frac{1}{\alpha_2 !} \leq \left( \frac{1}{2} \right)^{m - 2}.
\end{equation*}
Thus, we have
\begin{equation*}
\frac{1}{\alpha !} \leq \left( \frac{1}{2} \right)^{m - 2}.
\end{equation*}
Since $m \geq 3$ and $\displaystyle{\sum_{|\alpha | = m} 1 = m + 1}$, we find that
\begin{equation*}
\sum_{|\alpha | = m} \frac{1}{\alpha !} \leq ( m + 1)  \left( \frac{1}{2} \right)^{m - 2} \leq 2.
\end{equation*}
Therefore, we obtain \eqref{eq22}.

Let us now derive \eqref{eq21}. Fix $m \in \mathbb{N}$ and $f \in C^{ m + 1} ( \overline{\Omega}_H)$. Since $t y + (1 - t ) x \in \overline{\Omega}_H$ for $0 \leq t \leq 1$ and $x , y \in \Omega_H$, we apply \eqref{eq22} to check that
\begin{equation*}
| R_{ m + 1 }[f] (x,y) | \leq 2 (m + 1) | x - y |^{ m + 1 } | f |_{C^{m+1}}.
\end{equation*}
Therefore, the lemma follows.
\end{proof}

To derive basic properties of our weight function, we prepare the following lemma.
\begin{lemmaz}\label{lem22}
$(\mathrm{i})$ For each $i \in \{ 1, \cdots , N \}$ and $f \in L^1 ( B_{ r_\sigma } (a_i))$,
\begin{equation}\label{eq23}
\| f \|_{L^1 ( \sigma_i )} \leq \| f \|_{L^1 ( B_{ r_\sigma } (a_i) )} .
\end{equation}
$(\mathrm{ii})$ For each $x \in B_\delta ( a_k )$ and $g \in L^1 ( B_{\lambda h + r_\sigma } (x))$,
\begin{equation}\label{eq24}
\sum_{i \in \mathcal{R} (x , \lambda h)} \| g \|_{L^1 ( \sigma_i )} \leq \| g \|_{L^1 ( B_{\lambda h + r_\sigma } (x) \setminus \sigma_k )} .
\end{equation}
$(\mathrm{iii})$ For each $x \in B_\delta ( a_k ) $ and $0 < \lambda < 1$,
\begin{equation}\label{eq25}
\sum_{i \in \mathcal{R} (x , h)} \int_{\sigma_i \cap B_h (x)} \frac{| y - a_i |}{| x - a_i |} { \ } d y \leq \frac{ \pi r_\sigma h }{ \lambda } + \pi ( \lambda h + r_\sigma )^2 .
\end{equation}
\end{lemmaz}

\begin{proof}(Lemma \ref{lem22})
We first show $(\mathrm{i})$. By the definition of $r_\sigma$, we find that $y \in B_{r_\sigma} (a_i)$ if $y \in \sigma_i$. Therefore, we see that for each $i = 1, \cdots, N$ and $f \in L^1 (B_{r_\sigma} (a_i) )$,
\begin{equation*}
\int_{\sigma_i} f (y) { \ }d y \leq \int_{B_{r_\sigma} (a_i) } f (y) { \ } d y.
\end{equation*}
This is \eqref{eq23}. Next, we derive \eqref{eq24}. Let $x \in B_\delta ( a_k )$ and $g \in L^1 ( B_{\lambda h + r_\sigma } (x))$. By the definition of $\mathcal{R} ( x , \lambda h)$, we find that $a_i \in B_{\lambda h} (x)$ if $i \in \mathcal{R} (x , \lambda h )$. Therefore, we observe that
\begin{equation*}
\sum_{i \in \mathcal{R} ( x , \lambda h ) } \int_{\sigma_i} g (y) { \ }d y \leq \int_{B_{ \lambda h + r_\sigma } (x) \setminus \sigma_k } g (y) { \ }d y.
\end{equation*}
Finally, we prove $(\mathrm{iii})$. Let $i \in \mathcal{R} (x, h)$. Since $i \neq k$ and $x \in \sigma_k \cap \Omega$, we find that
\begin{equation}\label{eq26}
\sup_{y \in \sigma_i } \frac{| y - a_i |}{ | x - a_i |} \leq 1 \text{ if } i \in \mathcal{R} (x , h).
\end{equation}
Let $j \in \mathcal{R}(x,h) \setminus \mathcal{R} (x , \lambda h)$. By definition, we see that
\begin{equation*}
\lambda h \leq | x - a_j | < h. 
\end{equation*}
This shows that
\begin{equation}\label{eq27}
\frac{1}{ | x - a_j |} \leq \frac{1}{ \lambda h } \text{ if } j \in \mathcal{R}(x,h) \setminus \mathcal{R} (x , \lambda h).
\end{equation}
Since
\begin{equation*}
\sum_{i \in \mathcal{R} (x ,h)} F_i  = \sum_{j \in \mathcal{R} (x ,h) \setminus \mathcal{R} (x , \lambda h)} F_j + \sum_{i \in \mathcal{R} (x , \lambda h)} F_i,
\end{equation*}
we use \eqref{eq24} to check that
\begin{multline*}
\sum_{i \in \mathcal{R} (x , h)} \int_{\sigma_i \cap B_h (x)} \frac{| y - a_i |}{| x - a_i |} { \ } d y\\
  \leq \frac{ r_\sigma }{ \lambda h } \sum_{i \in \mathcal{R} (x,h) \setminus \mathcal{R} (x , \lambda h) } \int_{\sigma_i \cap B_h (x)} 1 { \ }d y + \sum_{i \in \mathcal{R} (x, \lambda h ) } \int_{\sigma_i} 1 { \ }d y \\
 \leq \frac{ \pi r_\sigma h }{ \lambda } + \pi ( \lambda h + r_\sigma )^2 .
\end{multline*}
Thus, we have \eqref{eq25}. Therefore, the lemma follows.
\end{proof}

\begin{lemmaz}[Properties of weight functions]\label{lem23}
Let $w$ be a weight function satisfying the properties of Assumption \ref{ass11}. Then the following six assertions hold:\\
$(\mathrm{i})$ For each $x \in B_\delta ( a_k )$ and $n \in \mathbb{Z}$,
\begin{equation}\label{eq28}
| x - \cdot |^n w ( x - \cdot ) \in L^1 ( B_{h, \delta } (x)  ).
\end{equation}
$(\mathrm{ii})$ For each $i= 1,2$, $x \in B_\delta ( a_k )$, and $n \in \mathbb{Z}$,
\begin{equation}\label{eq29}
\int_{ B_{h , \delta } (x) } ( x_i -y_i ) | x - y |^n w (x - y ) { \ } d y = 0.
\end{equation}
$(\mathrm{iii})$ For each $i , j = 1,2$, $x \in B_\delta ( a_k )$, and $n \in \mathbb{Z}$,
\begin{equation}\label{Eq210}
\int_{ B_{h, \delta } (x) } (x_i -y_i )( x_j - y_j ) | x - y |^n w (x - y ) { \ } d y = \frac{1}{2} \delta_{i j} \| | x - \cdot |^{n + 2} w ( x - \cdot ) \|_{L^1 (B_{h, \delta} (x) )}.
\end{equation}
Here $\delta_{i j}$ denotes the Kronecker delta.\\
$( \mathrm{iv})$ For each $x \in B_\delta ( a_k )$,
\begin{equation}\label{Eq211}
\left| \sum_{i \in \mathcal{R} (x , h)}  \int_{\sigma_i} w ( x - y) { \ } d y  - \sum_{i \in \mathcal{R} (x , h)} V_i (x) w ( x - a_i ) \right| \leq \pi L_w r_\sigma h^2.
\end{equation}
$(\mathrm{v})$ For each $x \in B_\delta ( a_k )$ and $0 < \lambda <1$,
\begin{equation}\label{Eq212}
\frac{ \displaystyle{\sum_{i \in \mathcal{R} (x ,h)} \int_{\sigma_i} \frac{| y - a_i|}{ | x - a_i| } w ( x - y ){ \ } d y} }{ \displaystyle{\sum_{j \in \mathcal{R} (x,h)}\int_{\sigma_j} w (x - y ) { \ } d y  } } \leq \frac{r_\sigma}{ \lambda h} + \frac{ \| w (x - \cdot ) \|_{L^1 (B_{\lambda h + r_\sigma } (x) \setminus \sigma_k )}}{ \| w ( x - \cdot ) \|_{L^1 ( B_h (x) \setminus \sigma_k )} }.
\end{equation}
$(\mathrm{vi})$ For each $x \in B_\delta ( a_k )$ and $0 < \lambda <1$,
\begin{multline}\label{Eq213}
\sum_{i \in \mathcal{R} (x ,h)} \int_{\sigma_i} \frac{| y - a_i|}{ | x - a_i| } \frac{w ( x - y )}{ | x - y| }{ \ } d y\\ \leq \left\| \frac{w (x - \cdot )}{| x - \cdot |} \right\|_{L^1 (B_{\lambda h + r_\sigma} (x) \setminus \sigma_k )} + \frac{r_\sigma}{ \lambda h} \left\| \frac{w ( x - \cdot )}{| x - \cdot |} \right\|_{L^1 ( B_h (x) \setminus \sigma_k )}.
\end{multline}
\end{lemmaz}

\begin{proof}(Lemma \ref{lem23})
We first show $(\mathrm{i})$. Using a change of variables, we see that
\begin{equation*}
\| w \|_{L^1 ( B_{h , \delta } )} = \int_{B_{h, \delta} } w (z) { \ }d z = \int_0^{2 \pi} \int_\delta^h \widehat{w} (r) r { \ }d r d \vartheta .
\end{equation*}
This gives
\begin{equation}\label{Eq214}
\int_\delta^h r \widehat{w} (r) { \ }d r = \frac{1}{2 \pi} \| w \|_{L^1 ( B_{h , \delta } ) }.
\end{equation}
Fix $x \in B_\delta (a_k)$ and $n \in \mathbb{Z}$. Using a change of variables with \eqref{Eq214}, we observe that
\begin{equation}\label{Eq215}
\| | x - \cdot |^n w (x - \cdot ) \|_{L^1 ( B_{h , \delta} (x) )} = 2 \pi \int_\delta^h r^{n + 1} \widehat{w} (r) { \ }d r \leq \max \{ h^n , \delta^n \} \| w \|_{L^1 (B_{h , \delta })} < + \infty .
\end{equation}
Therefore, we obtain $(\mathrm{i})$.

Next, we prove $( \mathrm{ii})$ and $(\mathrm{iii})$. Fix $i , j = 1,2$, $x \in B_\delta ( a_k )$, and $n \in \mathbb{Z}$. From \eqref{Eq215}, we have
\begin{equation}\label{Eq216}
\int_{\delta}^h r^{n+1} \widehat{w} (r) { \ }d r = \frac{1}{2 \pi} \| | x - \cdot |^n w (x - \cdot ) \|_{L^1 (B_{h , \delta } (x) )}.
\end{equation}
Using a change of variables, we find that
\begin{equation*}
\int_{B_{h , \delta} (x) } ( x_i - y_i) | x - y |^n w ( x - y ) { \ } d y = 0 \cdot \int_\delta^h r^{n+1} \widehat{w} (r) { \ }d r = 0.
\end{equation*}
By \eqref{Eq215}, we check that
\begin{multline*}
\int_{B_{h , \delta} (x) } ( x_i - y_i ) (x_j - y_j )| x - y |^n w ( x - y ) { \ } d y\\ = \int_{B_{h , \delta} (x) } \frac{ ( x_i - y_i ) (x_j - y_j )}{| x - y |^2}  | x - y |^{n + 2 } w ( x - y ) { \ } d y\\
= \pi \delta_{i j} \int_{\delta}^hr^{n + 1} \widehat{w} (r) { \ } d r = \frac{1}{2} \delta_{i j} \| | x - \cdot |^n w ( x - \cdot ) \|_{L^1 (B_{h , \delta} (x) )}.
\end{multline*}
Therefore, we have $(\mathrm{ii})$ and $(\mathrm{iii})$.

Now, we show $(\mathrm{iv})$. Since $\displaystyle{ B_h (x) \subset \sum_{i \in \overline{\mathcal{R}} (x ,h)} \overline{\sigma}_i }$ and $w (x ) =0$ for $x \in \mathbb{R}^2 \setminus \overline{B}_{h , \delta}$, we use the Lipschitz continuity of $w$ to observe that
\begin{align*}
\text{ (L.H.S.) of } \eqref{Eq211} & = \left| \sum_{i \in \mathcal{R} (x , h)}  \int_{\sigma_i \cap B_h (x)} \{ w ( x - y) - w (x - a_i ) \} { \ } d y \right|\\
& \leq L_w r_\sigma \sum_{i \in \mathcal{R} (x , h)}  \int_{\sigma_i \cap B_h (x)} 1 { \ } d y\\
& \leq \pi L_w r_\sigma h^2 .
\end{align*}
Thus, we have \eqref{Eq211}.

Finally, we prove $( \mathrm{v})$ and $( \mathrm{vi})$. Using \eqref{eq26}, \eqref{eq27}, and \eqref{eq24}, we check that
\begin{multline*}
\sum_{i \in \mathcal{R} (x ,h)} \int_{\sigma_i} \frac{| y - a_i|}{ | x - a_i| } w ( x - y ){ \ } d y\\ \leq \frac{r_\sigma}{ \lambda h} \| w ( x - \cdot ) \|_{L^1 ( B_h (x) \setminus \sigma_k )} + \| w (x - \cdot ) \|_{L^1 (B_{\lambda h + r_\sigma} (x) \setminus \sigma_k )}.
\end{multline*}
This gives \eqref{Eq212}. Similarly, we see \eqref{Eq213}. Therefore, the lemma follows.
\end{proof}

\section{Error Estimate $(\mathrm{I})$}\label{sect3}

In this section, we study $| f ( a_k ) - \widetilde{\Pi}_h f ( a_k ) |$ to prove Theorem \ref{thm13}. To this end, we introduce some notations. Let $w$ be a weight function satisfying the properties of Assumption \ref{ass11}. For each $f \in C^0 ( \overline{\Omega}_H )$, we define
\begin{align*}
\Pi_h f ( a_k ) & = \frac{ \displaystyle{ \int_{B_{h, \delta } ( a_k )} f (y) w ( a_k -y ) { \ }d y }  }{ \displaystyle{ \int_{B_{h , \delta } } w (z) { \ }d z }  } ,\\
\widehat{\Pi}_h f ( a_k ) & = \frac{ \displaystyle{ \sum_{i \in \mathcal{R} ( a_k ,h) }\int_{\sigma_i} f (y) w ( a_k - y) { \ }d y }  }{ \displaystyle{ \sum_{j \in \mathcal{R} ( a_k , h) }\int_{\sigma_j} w ( a_k - z) { \ }d z }  },\\
\breve{\Pi}_h f ( a_k ) & = \frac{ \displaystyle{ \sum_{i \in \mathcal{R} ( a_k , h) } V_i ( a_k ) f (a_i) w ( a_k  - a_i ) }  }{ \displaystyle{ \sum_{j \in \mathcal{R} ( a_k , h) } \int_{\sigma_j} w ( a_k  - z ) { \ } d z }  } .
\end{align*}
It is easy to check that $| \Pi_h f ( a_k ) | \leq | f |_{C^0}$. Since
\begin{multline*}
f ( a_k ) - \widetilde{\Pi}_h f ( a_k ) = \{ f ( a_k ) - \Pi_h f ( a_k ) \} + \{ \Pi_h f ( a_k ) - \widehat{\Pi}_h f ( a_k ) \}\\
 + \{ \widehat{\Pi}_h f ( a_k ) - \breve{\Pi}_h f ( a_k ) \} + \{ \breve{\Pi}_h f ( a_k ) - \widetilde{\Pi}_h f ( a_k ) \},
\end{multline*}
we prove the following lemma.
\begin{lemmaz}\label{lem31}
For each $f \in C^1 (\overline{\Omega}_H)$,
\begin{align}
| f ( a_k ) - \Pi_h f ( a_k ) | & \leq h | f |_{C^1} ,\label{eq31}\\
| \Pi_h f ( a_k ) - \widehat{\Pi}_h f ( a_k ) | & \leq 2 \frac{ \| w ( a_k  - \cdot ) \|_{L^1 (\sigma_k \setminus B_{\delta} (a_k) )}}{ \| w (  a_k  - \cdot ) \|_{L^1 (B_h ( a_k ) \setminus B_\delta (a_k) )}} | f |_{C^0} , \label{eq32}\\
| \widehat{\Pi}_h f ( a_k ) - \breve{\Pi}_h f ( a_k ) | & \leq r_\sigma | f |_{C^1} + \frac{\pi L_w r_\sigma h^2}{ \displaystyle{ \| w (  a_k  - \cdot )  \|_{L^1 (B_h ( a_k ) \setminus \sigma_k)} } } | f |_{C^0}, \label{eq33}\\
| \breve{\Pi}_h f ( a_k ) - \widetilde{\Pi}_h f ( a_k ) | & \leq \frac{ \pi L_w r_\sigma h^2}{ \displaystyle{ \| w (  a_k  - \cdot ) \|_{L^1 (B_h ( a_k ) \setminus \sigma_k )} }  } | f |_{C^0}. \label{eq34}
\end{align}
\end{lemmaz}

\begin{proof}(Lemma \ref{lem31})
We first show \eqref{eq31}. Since
\begin{equation*}
f ( a_k ) = \frac{1}{ \| w \|_{L^1 ( B_{h, \delta })}} \int_{B_{h, \delta } ( a_k )} f (  a_k  ) w (  a_k -y ) { \ }d y,
\end{equation*}
we use the mean-value theorem to see that
\begin{align*}
| f ( a_k ) - \Pi_h f ( a_k ) | & \leq \frac{ 1 }{ \| w \|_{L^1 ( B_{h , \delta } ) }} \int_{B_{h, \delta } ( a_k )} |f ( a_k ) - f (y)| w ( a_k -y) { \ }d y\\
& \leq h | f |_{C^1}.
\end{align*}
Therefore, we have \eqref{eq31}.

Secondly, we derive \eqref{eq32}. From 
\begin{equation*}
\int_{B_{h , \delta} (a_k) } w ( a_k - y ) { \ }d y = \sum_{i \in \overline{\mathcal{R}} (a_k , h) } \int_{\sigma_i \cap B_{h , \delta} (a_k) } w ( a_k - y ) { \ }d y,
\end{equation*}
we observe that
\begin{equation*}
| \Pi_h f ( a_k ) - \widehat{\Pi}_h f ( a_k ) | \leq 2 \frac{ | f |_{C^0} \| w ( a_k  - \cdot ) \|_{L^1 (\sigma_k \setminus B_{\delta} (a_k) )}}{ \| w (  a_k  - \cdot ) \|_{L^1 (B_h ( a_k ) \setminus B_\delta (a_k) )}} .
\end{equation*}
Therefore, we see \eqref{eq32}.

Next, we prove \eqref{eq33}. Since
\begin{multline*}
\int_{\sigma_i} f (y) w (  a_k  - y) { \ } d y - V_i ( a_k ) f (a_i) w (  a_k  - a_i ) \\
= \int_{\sigma_i \cap B_{h , \delta } ( a_k )} \{ f (y) - f (a_i)\} w (  a_k  - y) { \ } d y\\ + \int_{\sigma_i \cap B_{h , \delta } ( a_k )} f (a_i) \{ w (  a_k  - y) - w (  a_k  - a_i ) \} { \ }d y,
\end{multline*}
we use the mean-value theorem and \eqref{Eq211} to check that
\begin{equation*}
| \widehat{\Pi}_h f ( a_k ) - \breve{\Pi}_h f ( a_k ) | \leq r_\sigma | f |_{C^1} + \frac{\pi L_w r_\sigma h^2}{ \displaystyle{ \| w (  a_k  - \cdot )  \|_{L^1 (B_h ( a_k ) \setminus \sigma_k)} } } | f |_{C^0}. 
\end{equation*}
This is \eqref{eq33}.

Finally, we show \eqref{eq34}. By \eqref{Eq211}, we find that
\begin{equation*}
| \breve{\Pi}_h f ( a_k ) - \widetilde{\Pi}_h f ( a_k  )| \leq \frac{ \pi L_w r_\sigma h^2 }{ \displaystyle{ \| w (  a_k  - \cdot ) \|_{L^1 ( B_h ( a_k ) \setminus \sigma_k )}  } } | f |_{C^0}.
\end{equation*}
Therefore, we have \eqref{eq34}, and the lemma follows.
\end{proof}

Finally, we prove Theorem \ref{thm13}.
\begin{proof}(Theorem \ref{thm13})
Using Lemma \ref{lem31}, we prove Theorem \ref{thm13}.
\end{proof}

\section{Error Estimate $(\mathrm{II})$}\label{sect4}

In this section, we consider $| \nabla f (a_k ) - \widetilde{\nabla}_h f (a_k) |$ to prove Theorem \ref{thm14}. Let $w$ be a weight function satisfying the properties of Assumption \ref{ass11}. For each $f \in C^0 ( \overline{\Omega}_H )$, we define
\begin{align*}
\nabla_h f ( a_k ) & = 2 \frac{ \displaystyle{ \int_{B_{h , \delta } ( a_k )} \frac{ f ( a_k ) - f (y) }{ |  a_k  - y| } \frac{ a_k -y}{ |  a_k  - y |} w ( a_k -y) { \ }d y }  }{ \displaystyle{ \int_{B_{h , \delta }} w (z) { \ }d z }  } ,\\
\widehat{\nabla}_h f ( a_k ) & = 2 \frac{ \displaystyle{ \sum_{i \in \mathcal{R} ( a_k ,h) }\int_{\sigma_i} \frac{ f ( a_k ) - f (y) }{ |  a_k  - y| } \frac{ a_k -y}{ |  a_k  - y |}  w ( a_k -y) { \ }d y }  }{ \displaystyle{ \sum_{j \in \mathcal{R} ( a_k ,h) } \int_{\sigma_j} w ( a_k - z ) { \ }d z }  },\\
\breve{\nabla}_h f ( a_k ) & = 2 \frac{ \displaystyle{ \sum_{i \in \mathcal{R} ( a_k ,h) }\int_{\sigma_i} \frac{ f ( a_k ) - f ( a_i ) }{ |  a_k  - a_i | } \frac{ a_k - a_i}{ |  a_k  - a_i |} w ( a_k - y) { \ }d y }  }{ \displaystyle{ \sum_{j \in \mathcal{R} ( a_k ,h) } \int_{\sigma_j}  w (  a_k  - z ) { \ }d z }  } .
\end{align*}
It is easy to check that
\begin{equation*}
\max  \{ | \nabla_h f ( a_k ) ,  | \widehat{\nabla}_h f ( a_k ) | , | \breve{\nabla}_h f ( a_k ) | \} | \leq \frac{4}{\delta} | f |_{C^0}.
\end{equation*}
It is clear that
\begin{multline*}
\nabla f ( a_k ) - \widetilde{\nabla}_h f ( a_k ) = \{ \nabla f ( a_k ) - \nabla_h f ( a_k ) \} + \{ \nabla_h f ( a_k ) - \widehat{\nabla}_h f ( a_k ) \}\\ + \{ \widehat{\nabla}_h f ( a_k ) - \breve{\nabla}_h f ( a_k ) \} + \{ \breve{\nabla}_h f ( a_k ) - \widetilde{\nabla}_h f ( a_k ) \} .
\end{multline*}

The aim of this section is to prove the following two lemmas.
\begin{lemmaz}\label{lem41}For each $f \in C^2 ( \overline{\Omega}_H)$,
\begin{equation}\label{eq41}
|\nabla f ( a_k ) - \nabla_h f ( a_k ) | \leq 4 h | f |_{C^2}.
\end{equation}
\end{lemmaz}

\begin{lemmaz}\label{lem42}For each $f \in C^1 ( \overline{\Omega}_H)$,
\begin{align}
| \nabla_h f ( a_k ) - \widehat{\nabla}_h f ( a_k ) | & \leq 4 \frac{ \| w (  a_k  - \cdot ) \|_{L^1 ( \sigma_k \setminus B_\delta (a_k) )}}{\| w ( a_k  - \cdot ) \|_{L^1 ( B_h ( a_k ) \setminus B_\delta (a_k) )}} | f |_{C^1} , \label{eq42}\\
| \widehat{\nabla}_h f ( a_k ) - \breve{\nabla}_h f ( a_k ) | & \leq 8 \left( \frac{r_\sigma}{ \lambda h} + \frac{ \displaystyle{ \| w (  a_k  - \cdot ) \|_{L^1 ( B_{ \lambda h + r_\sigma } ( a_k ) \setminus B_{\delta} ( a_k )  )} }  }{ \displaystyle{ \| w (  a_k  - \cdot ) \|_{L^1 (B_h ( a_k ) \setminus \sigma_k )} }  }  \right) | f |_{C^1}, \label{eq43}\\
| \breve{\nabla}_h f ( a_k ) - \widetilde{\nabla}_h f ( a_k ) | & \leq 4 \frac{ \pi L_w  r_\sigma h^2 }{ \| w (a_k - \cdot ) \|_{L^1 ( B_h (a_k) \setminus \sigma_k )} } | f |_{C^1}. \label{eq44}
\end{align}
\end{lemmaz}

We first show Lemma \ref{lem41}. Then we prove Lemma \ref{lem42}.
\begin{proof}(Lemma \ref{lem41})
Fix  $ f \in C^2 ( \overline{\Omega}_H)$. From the Taylor expansion, we have
\begin{equation}\label{eq45}
\sum_{| \alpha| =1} \frac{D^\alpha f ( a_k )}{\alpha !} (  a_k  - y )^\alpha - \{ f ( a_k ) - f (y) \} = R_2 [ f ] ( a_k  , y ),
\end{equation}
where $y \in B_{h , \delta} (a_k)$. Here
\begin{equation*}
R_2 [ f ] ( a_k ,y) = \sum_{| \alpha | = 2} (y -  a_k  )^\alpha \frac{2}{\alpha !} \int_0^1 D^\alpha f ( t y + (1 - t )  a_k  ) { \ } d t.
\end{equation*}
Multiplying both sides of \eqref{eq45} by $2 w ( a_k -y)( a_k -y)/( |  a_k  -y |^2 \| w \|_{L^1 (B_{h , \delta })} )$, and then integrating with respect to $y$, we have
\begin{equation}\label{eq46}
P_1 ( a_k ) - \widetilde{\nabla}_h f ( a_k ) = P_2 ( a_k ).
\end{equation}
Here
\begin{align*}
P_1 ( a_k ) &:= \frac{2}{\| w \|_{L^1 (B_{h, \delta })}} \int_{B_{h, \delta } ( a_k )} \sum_{| \alpha| =1} \frac{D^\alpha f ( a_k )}{\alpha !} (  a_k  - y )^\alpha \frac{( a_k -y)}{| a_k -y|^2} w ( a_k -y) { \ }d y,\\
P_2 ( a_k ) &:= \frac{2}{\| w \|_{L^1 (B_{h , \delta })}}  \int_{B_{h , \delta } ( a_k )} R_2 [f] ( a_k  , y ) \frac{( a_k -y)}{| a_k -y|^2} w ( a_k -y) { \ }d y.
\end{align*}
From \eqref{Eq210}, we see that
\begin{equation}\label{eq47}
P_1 ( a_k ) = \nabla f ( a_k ) .
\end{equation}
By \eqref{eq21}, we find that
\begin{equation}\label{eq48}
| P_2 ( a_k ) | \leq 4 h | v |_{C^2} .
\end{equation}
Combining \eqref{eq46}-\eqref{eq48} gives \eqref{eq41}. Therefore, the lemma follows.
\end{proof}

\begin{proof}(Lemma \ref{lem42})
Let $f \in C^1 (\overline{\Omega}_H)$. We first show \eqref{eq42}. By the mean-value theorem, we see that
\begin{multline*}
\frac{1}{2} | \nabla_h f ( a_k ) - \widehat{\nabla}_h f ( a_k ) |\\
\leq \frac{1}{ \| w (  a_k  - \cdot ) \|_{L^1 ( B_h ( a_k ) \setminus B_\delta (a_k) ) }} \left| \int_{\sigma_k} \frac{ f ( a_k ) - f (y) }{ |  a_k  - y| } \frac{ a_k -y}{ |  a_k  - y |} w ( a_k -y) { \ }d y \right|\\
+ \frac{ \| w (  a_k  - \cdot )\|_{L^1 ( \sigma_k )}  }{ \| w ( a_k  - \cdot ) \|_{L^1 ( B_h ( a_k ) \setminus \sigma_k )} \| w \|_{L^1 (B_{h , \delta }) }} \cdot\\
\left| \sum_{i \in \mathcal{R} (a_k , \lambda )} \int_{\sigma_i} \frac{ f ( a_k ) - f (y) }{ |  a_k  - y| } \frac{ a_k -y}{ |  a_k  - y |} w ( a_k -y) { \ }d y \right|\\
\leq 2 \frac{ \| w (  a_k  - \cdot ) \|_{L^1 ( \sigma_k \setminus B_\delta (a_k) )}}{\| w ( a_k  - \cdot ) \|_{L^1 ( B_h ( a_k ) \setminus B_\delta (a_k) )} } | f |_{C^1}.
\end{multline*}
Therefore, we have \eqref{eq42}. Note that $\| w ( a_k  - \cdot ) \|_{L^1 ( \sigma_k )} = \| w ( a_k  - \cdot ) \|_{L^1 ( \sigma_k \setminus B_\delta (a_k) )}$.

Next, we derive \eqref{eq43}. A direct calculation gives
\begin{equation}\label{eq49}
\frac{1}{2} \{ \widehat{\nabla}_h f ( a_k ) - \breve{\nabla}_h f ( a_k ) \} = P_3 ( a_k ) + P_4 ( a_k ) + P_5 ( a_k ) + P_6 ( a_k ).
\end{equation}
Here
\begin{multline*}
P_3 ( a_k ) : =\\ \frac{ \displaystyle{ \sum_{i \in \mathcal{R} ( a_k ,h) }\int_{\sigma_i} \{ f ( a_k ) - f (y) \} \frac{( a_k -y)}{ |  a_k  - y | } \left( \frac{1}{ |  a_k  - y | } - \frac{1}{ |  a_k  - a_i | } \right) w ( a_k -y) { \ }d y }  }{ \displaystyle{ \sum_{j \in \mathcal{R} ( a_k ,h) }\int_{\sigma_j} w ( a_k - z ) { \ }d z }  },
\end{multline*}
\begin{equation*}
P_4 ( a_k ) := \frac{ \displaystyle{ \sum_{i \in \mathcal{R} ( a_k ,h) }\int_{\sigma_i} \{ f (a_i) - f (y) \} \frac{( a_k -y)}{ |  a_k  - y | |  a_k  - a_i| } w ( a_k -y) { \ }d y }  }{ \displaystyle{ \sum_{j \in \mathcal{R} ( a_k ,h) }\int_{\sigma_j} w ( a_k - z ) { \ }d z }  } ,
\end{equation*}
\begin{multline*}
P_5 ( a_k ) :=\\ \frac{ \displaystyle{ \sum_{i \in \mathcal{R} ( a_k ,h) }\int_{\sigma_i} \frac{ f ( a_k ) - f (a_i) }{ |  a_k  - a_i| } ( a_k -y) \left( \frac{1}{ |  a_k  - y | } - \frac{1}{ |  a_k  - a_i | } \right) w ( a_k -y) { \ }d y }  }{ \displaystyle{ \sum_{j \in \mathcal{R} ( a_k ,h) }\int_{\sigma_j} w ( a_k - z ) { \ }d z }  } ,
\end{multline*}
and
\begin{equation*}
P_6 ( a_k ) := \frac{ \displaystyle{ \sum_{i \in \mathcal{R} ( a_k ,h) }\int_{\sigma_i} \frac{ f ( a_k ) - f (a_i) }{ |  a_k  - a_i|^2 } (a_i -y) w ( a_k -y) { \ }d y }  }{ \displaystyle{ \sum_{j \in \mathcal{R} ( a_k ,h) }\int_{\sigma_j} w ( a_k - z ) { \ }d z }  } .
\end{equation*}
Since
\begin{equation*}
\left| \frac{1}{ |  a_k  - y|} - \frac{1}{ |  a_k - a_i |} \right| \leq \frac{ | a_i - y |}{ |  a_k  - y| |  a_k  - a_i| }, 
\end{equation*}
we use the mean-value theorem and \eqref{Eq212} to check that
\begin{align}
\sum_{q=3}^6 P_q ( a_k ) & \leq 4| f |_{C^1}  \frac{ \displaystyle{ \sum_{i \in \mathcal{R} ( a_k ,h) }\int_{\sigma_i} \frac{ | y - a_i| }{ |  a_k  - a_i |} w ( a_k -y) { \ }d y }  }{ \displaystyle{ \sum_{j \in \mathcal{R} ( a_k ,h) }\int_{\sigma_j} w ( a_k - z ) { \ }d z }  }\notag\\
& \leq 4 \left( \frac{r_\sigma}{\lambda h} +\frac{ \| w ( a_k - \cdot ) \|_{L^1 ( B_{\lambda h + r_\sigma } (a_k) \setminus \sigma_k )}  }{ \displaystyle{ \| w (  a_k  - \cdot ) \|_{L^1 (B_h ( a_k ) \setminus \sigma_k )} }  } \right) | f |_{C^1}.\label{Eq410}
\end{align}
From \eqref{eq49} and \eqref{Eq410}, we have \eqref{eq43}.

Finally, we show \eqref{eq44}. From \eqref{Eq211}, we see that
\begin{align*}
\frac{1}{2} | \breve{\nabla}_h f ( a_k ) - \widetilde{\nabla}_h f ( a_k ) | \leq 2 \frac{ \pi L_w  r_\sigma h^2  | f |_{C^1} }{  \| w (a_k - \cdot ) \|_{L^1 (B_h (a_k) \setminus \sigma_k )} } .
\end{align*}
This is \eqref{eq44}. Therefore, the lemma follows.
\end{proof}

Finally, we prove Theorem \ref{thm14}.
\begin{proof}(Theorem \ref{thm14})
Combining Lemmas \ref{lem41} and \ref{lem42} gives Theorem \ref{thm14}.
\end{proof}

\section{Error Estimate $(\mathrm{III})$}\label{sect5}
In this section, we estimate $| \Delta f ( a_k ) - \widetilde{\Delta }_h f (a_k) |$ to prove Theorem \ref{thm15}. Let $w$ be a weight function satisfying the properties of Assumption \ref{ass11}. For each $f \in C^0 ( \overline{\Omega}_H )$, we define
\begin{align*}
\Delta_h f (a_k) & = - 4 \frac{ \displaystyle{ \int_{B_{h, \delta } (a_k)} \{ f (a_k) - f (y) \} w (a_k-y) { \ }d y }  }{ \displaystyle{ \int_{B_{h , \delta }} | z |^2 w (z) { \ }d z }  } ,\\
\widehat{\Delta}_h f (a_k) & = - 4 \frac{ \displaystyle{ \sum_{i \in \mathcal{R} (a_k,h) }\int_{\sigma_i} \{ f (a_k) - f (y) \} w (a_k - y ) { \ }d y }  }{ \displaystyle{ \sum_{j \in \mathcal{R} (a_k,h) }\int_{\sigma_j} | a_k - z |^2 w (a_k - z ) { \ }d z }  },\\
\breve{\Delta}_h f (a_k) & = - 4 \frac{ \displaystyle{ \sum_{i \in \mathcal{R} (a_k,h) }\int_{\sigma_i} \{ f (a_k) - f ( a_i ) \} w (a_k- y ) { \ }d y }  }{ \displaystyle{ \sum_{j \in \mathcal{R} (a_k,h) } \int_{\sigma_j} | a_k - z |^2  w ( a_k - z ) { \ }d z }  } .
\end{align*}
We see at once that
\begin{equation*}
| \Delta_h f (a_k)| \leq \frac{ 8 }{\delta^2 } | f |_{C^0} < + \infty .
\end{equation*}

The proof of Theorem \ref{thm15} makes use of the following two lemmas.
\begin{lemmaz}\label{lem51}For each $f \in C^3 ( \overline{\Omega}_H)$,
\begin{equation}\label{eq51}
| \Delta f (a_k) - \Delta_h f (a_k) | \leq 24 h | f |_{C^3}.
\end{equation}
\end{lemmaz}

\begin{lemmaz}\label{lem52}For each $f \in C^1 ( \overline{\Omega}_H)$,
\begin{align}
| \Delta_h f (a_k) - \widehat{\Delta}_h f (a_k) | & \leq 4 c_4 (a_k) | f |_{C^1},\label{eq52}\\
| \widehat{\Delta}_h f (a_k) - \breve{\Delta}_h f (a_k) | & \leq \{ 4 c_5 (a_k) + 4 c_6 (a_k) \} | f |_{C^1}, \label{eq53}\\
| \breve{\Delta}_h f (a_k) - \widetilde{\Delta}_h f (a_k) | & \leq \{ 4 c_7 (a_k) + 4 c_8 (a_k) \} | f |_{C^1}. \label{eq54}
\end{align}
Here $c_4 (a_k), \cdots, c_8 (a_k)$ are the constants defined by Theorem \ref{thm15}.
\end{lemmaz}

We first show Lemma \ref{lem51}. Then we prove Lemma \ref{lem52}.
\begin{proof}(Lemma \ref{lem51})
Fix  $f \in C^3 ( \overline{\Omega}_H )$. From the Taylor expansion, we have
\begin{equation}\label{eq55}
- \sum_{| \alpha| = 2} \frac{D^\alpha f ( a_k )}{\alpha !} (  a_k  - y )^\alpha - \{ f ( a_k ) - f (y) \} = - \sum_{| \alpha| =1} \frac{D^\alpha f ( a_k )}{\alpha !} (  a_k  - y )^\alpha + R_3 [f] ( a_k  , y ),
\end{equation}
where $y \in B_{h , \delta} ( a_k ) $. Here
\begin{equation*}
R_{ 3 }[f] ( a_k ,y) := \sum_{| \alpha | = 3} (y -  a_k  )^\alpha \frac{3}{\alpha !} \int_0^1 (1 - t)^2 D^\alpha f ( t y + (1 - t )  a_k  ) { \ } d t.
\end{equation*}
Multiplying both sides of \eqref{eq55} by $ -4 w (a_k - y) / \| | \cdot |^2 w (\cdot ) \|_{L^1 ( B_{h , \delta} )}$, and integrating with respect to $y$, we have
\begin{equation}\label{eq56}
Q_1 (a_k) - \Delta_h f (a_k) = Q_2 (a_k) + Q_3 (a_k) .
\end{equation}
Here
\begin{align*}
Q_1 (a_k) &:= \frac{4}{ \| | \cdot |^2 w( \cdot ) \|_{L^1 (B_{h, \delta })}} \int_{B_{h, \delta } (a_k)} \sum_{| \alpha| =2} \frac{D^\alpha f (a_k)}{\alpha !} ( a_k - y )^\alpha w (a_k-y) { \ }d y,\\
Q_2 (a_k) &:= \frac{ 4 }{ \| | \cdot |^2 w ( \cdot ) \|_{L^1 (B_{h, \delta })}} \int_{B_{h, \delta } (a_k)} \sum_{| \alpha| =1} \frac{D^\alpha f (a_k)}{\alpha !} ( a_k - y )^\alpha w (a_k-y) { \ }d y,\\
Q_3 (a_k) &:= - \frac{4}{ \| | \cdot |^2 w (\cdot ) \|_{L^1 (B_{h , \delta }) }} \int_{B_{h, \delta} (a_k)} R_3 [f] (a_k , y ) w (a_k-y) { \ }d y.
\end{align*}
Using \eqref{Eq210}, \eqref{eq29}, and \eqref{eq21}, we find that
\begin{align}
Q_1 (a_k) & = \Delta f (a_k),\label{eq57}\\
Q_2 (a_k) & = 0 , \label{eq58}\\
| Q_3 (a_k) | & \leq 24 h | f |_{C^3} .\label{eq59}
\end{align}
Combining \eqref{eq56}-\eqref{eq59}, we have \eqref{eq51}. Therefore, the lemma follows.
\end{proof}

\begin{proof}(Lemma \ref{lem52})
Let $f \in C^1 ( \overline{\Omega}_H)$. We first show \eqref{eq52}. Since
\begin{equation*}
| a_k - y | \leq r_\sigma 
\end{equation*}
for $y \in \sigma_k$, we use the mean-value theorem to see that
\begin{multline*}
\frac{1}{4} | \Delta_h f (a_k) - \widehat{\Delta}_h f (a_k) | \leq \frac{ \| | a_k - \cdot | w (a_k - \cdot ) \|_{L^1 (\sigma_k \setminus B_\delta (a_k) ) } }{ \| | a_k - \cdot |^2 w (a_k - \cdot ) \|_{L^1 (B_h (a_k) \setminus B_\delta (a_k) )} } \cdot \\
\bigg( 1 + r_\sigma \frac{ \| | a_k - \cdot | w (a_k - \cdot ) \|_{L^1 ( B_h (a_k) \setminus B_\delta (a_k) ) } }{ \| | a_k - \cdot |^2 w (a_k - \cdot ) \|_{L^1 (B_h (a_k) \setminus \sigma_k )} } \bigg) | f |_{C^1}.
\end{multline*}

Using $\sum_{i \in \mathcal{R} (a_k, h) } = \sum_{i \in \mathcal{R} (a_k , \lambda h)} + \sum_{i \in \mathcal{R} (a_k , h) \setminus \mathcal{R} (a_k , \lambda h) }$, we see that
\begin{multline*}
\frac{1}{4} | \widehat{\Delta}_h f (a_k) - \breve{\Delta}_h f (a_k) | \leq | f |_{C^1} \frac{ \displaystyle{\sum_{ i \in \mathcal{R} (a_k ,h ) }  \int_{\sigma_i } \frac{| a_i - y | }{ |a_k -y| } | a_k - y | w ( a_k - y ) { \ } d y } }{ \displaystyle{ \| | a_k - \cdot |^2 w ( a_k - \cdot ) \|_{L^1 (B_h (a_k) \setminus \sigma_k ) } }} \\
\leq \frac{ r_\sigma }{ \lambda h} \frac{ \displaystyle{ \| | a_k - \cdot | w (a_k - \cdot ) \|_{L^1 ( B_h (a_k) \setminus \sigma_k )} }  }{ \displaystyle{ \| | a_k - \cdot |^2 w ( a_k - \cdot ) \|_{L^1 (B_h (a_k) \setminus \sigma_k ) } }  } | f |_{C^1} 
\\ + r_\sigma \frac{ \displaystyle{ \| w (a_k - \cdot ) \|_{L^1 ( B_{\lambda h + r_\sigma} (a_k) \setminus \sigma_k )} }  }{ \displaystyle{ \| | a_k - \cdot |^2 w ( a_k - \cdot ) \|_{L^1 (B_h (a_k) \setminus \sigma_k ) } }  } | f |_{C^1} .
\end{multline*}
A direct calculation shows that
\begin{multline}\label{Eq510}
\bigg| |a_k - y|^2 w (a_k - y) - | a_k - a_i |^2 w (a_k - a_i) \bigg|\\
\leq \{ | a_k - y| + | a_k - a_i | \} | a_i - y | w (a_k -y) + | a_k - a_i|^2 \{ w (a_k - y) - w ( a_k - a_i )\}\\
\leq 2 r_\sigma h + L_w r_\sigma h^2.
\end{multline}
Using \eqref{Eq510}, we check that
\begin{multline*}
\frac{1}{4} | \breve{\Delta}_h f (a_k) - \widetilde{\Delta}_h f (a_k) |
\leq \frac{ \pi L_w r_\sigma h^3 }{ \| | a_k - \cdot |^2 w ( a_k - \cdot ) \|_{L^1 (B_h (a_k) \setminus \sigma_k )} } | f |_{C^1}\\
+ \frac{2 r_\sigma h \| w (a_k - \cdot ) \|_{L^1 (B_h (a_k) \setminus \sigma_k )} + \pi L_w r_\sigma h^4 }{ \| | a_k - \cdot |^2 w (a_k - \cdot ) \|_{L^1 (B_h (a_k) \setminus \sigma_k )}} \frac{ \displaystyle{\sum_{i \in \mathcal{R} (a_k , h )} V_i | a_k - a_i | w (a_k - a_i)}}{ \displaystyle{\sum_{j \in \mathcal{R} (a_k , h )} V_j | a_k - a_j |^2 w (a_k - a_j)  } } | f |_{C^1}.
\end{multline*}
Therefore, the lemma follows.
\end{proof}

Finally, we prove Theorem \ref{thm15}.
\begin{proof}(Theorem \ref{thm15})
Using Lemmas \ref{lem51} and \ref{lem52}, we prove Theorem \ref{thm15}.
\end{proof}

\section{Error Estimate $(\mathrm{IV})$}\label{sect6}
In this section, we consider $| \Delta f (a_k ) -\widetilde{ \square}_h f ( a_k ) |$ to prove Theorem \ref{thm16}. Let $w$ be a weight function satisfying the properties of Assumption \ref{ass11}. For each $f \in C^0 ( \overline{\Omega}_H )$, we define
\begin{align*}
\square_h f ( a_k ) & = - 4 \frac{ \displaystyle{ \int_{B_{h, \delta } ( a_k )} \frac{ f ( a_k ) - f (y) }{ |  a_k  - y|^2 } w ( a_k -y) { \ }d y }  }{ \displaystyle{ \int_{B_{h , \delta }} w (z) { \ }d z }  } ,\\
\widehat{\square}_h f ( a_k ) & = - 4 \frac{ \displaystyle{ \sum_{i \in \mathcal{R} ( a_k ,h) }\int_{\sigma_i} \frac{ f ( a_k ) - f (y) }{ |  a_k  - y|^2 }  w ( a_k -y) { \ }d y }  }{ \displaystyle{ \sum_{j \in \mathcal{R} ( a_k ,h) }\int_{\sigma_j} w ( a_k - z ) { \ }d z }  },\\
\breve{\square}_h f ( a_k ) & = - 4 \frac{ \displaystyle{ \sum_{i \in \mathcal{R} ( a_k ,h) }\int_{\sigma_i} \frac{ f ( a_k ) - f ( a_i ) }{ |  a_k  - a_i |^2 } w ( a_k - y) { \ }d y }  }{ \displaystyle{ \sum_{j \in \mathcal{R} ( a_k ,h) } \int_{\sigma_j}  w (  a_k  - z) { \ }d z }  } .
\end{align*}
It is easy to check that
\begin{equation*}
| \square_h f ( a_k )| \leq \frac{ 8 }{ \delta^2 } | f |_{C^0} < + \infty .
\end{equation*}

Let us attack the following two lemmas.
\begin{lemmaz}\label{lem61}For each $f \in C^3 ( \overline{\Omega}_H)$,
\begin{equation}\label{eq61}
| \Delta f ( a_k ) - \square_h f ( a_k ) | \leq 24 h | f |_{C^3}.
\end{equation}
\end{lemmaz}

\begin{lemmaz}\label{lem62}For each $f \in C^1 ( \overline{\Omega}_H)$,
\begin{align}
| \square_h f ( a_k ) - \widehat{\square}_h f ( a_k ) | & \leq 8 c_9 ( a_k ) | f |_{C^1}, \label{eq62}\\
| \widehat{\square}_h f ( a_k ) - \breve{\square}_h f ( a_k ) | & \leq 16 \{ c_{10} ( a_k ) + c_{11} (a_k) \} | f |_{C^1},\label{eq63}\\
| \breve{\square}_h f ( a_k ) - \widetilde{\square}_h f ( a_k ) | & \leq 4 c_{12} ( a_k ) | f |_{C^1}. \label{eq64}
\end{align}
Here $c_9 ( a_k ), \cdots, c_{12} ( a_k )$ are the constants defined by Theorem \ref{thm16}.
\end{lemmaz}

We first show Lemma \ref{lem61}. Then we prove Lemma \ref{lem62}.
\begin{proof}(Lemma \ref{lem61})
Fix  $f \in C^3 ( \overline{\Omega}_H )$. Multiplying both sides of \eqref{eq55} by $ - 4 w ( a_k  - y)/(| a_k -y|^2 \| w \|_{L^1 ( B_{h , \delta} )} )$, and then integrating with respect to $y$, we have
\begin{equation}\label{eq65}
O_1 ( a_k ) - \square_h f ( a_k ) = O_2 ( a_k ) + O_3 ( a_k ) .
\end{equation}
Here
\begin{align*}
O_1 ( a_k ) &:= \frac{4}{ \| w \|_{L^1 (B_{h, \delta })}} \int_{B_{h, \delta } ( a_k )} \sum_{| \alpha| =2} \frac{D^\alpha f ( a_k )}{\alpha !}  \frac{ (  a_k  - y )^\alpha }{| a_k -y|^2} w ( a_k -y) { \ }d y,\\
O_2 ( a_k ) &:= \frac{ 4 }{ \| w \|_{L^1 (B_{h, \delta })}} \int_{B_{h, \delta } ( a_k )} \sum_{| \alpha| =1} \frac{D^\alpha f ( a_k )}{\alpha !} \frac{ (  a_k  - y )^\alpha }{| a_k -y|^2} w ( a_k -y) { \ }d y,\\
O_3 ( a_k ) &:= - \frac{4}{ \| w \|_{L^1 (B_{h , \delta }) }} \int_{B_{h, \delta} ( a_k )} R_3 [f] ( a_k  , y ) \frac{1}{| a_k -y|^2} w ( a_k -y) { \ }d y.
\end{align*}
By \eqref{Eq210} and \eqref{eq29}, we find that
\begin{align}
O_1 ( a_k ) & = \Delta f ( a_k ),\label{eq66}\\
O_2 ( a_k ) & = 0.\label{eq67}
\end{align}
Applying \eqref{eq21}, we see that
\begin{equation}\label{eq68}
| O_3 ( a_k ) | \leq 24 h | f |_{C^3} .
\end{equation}
Combining \eqref{eq65}-\eqref{eq68}, we have \eqref{eq61}. Therefore, the lemma follows.
\end{proof}

\begin{proof}(Lemma \ref{lem62})
Let $f \in C^1 ( \overline{\Omega}_H)$. We first show \eqref{eq62}. By the mean-value theorem, we see that
\begin{multline*}
\frac{1}{4} | \square_h f ( a_k ) - \widehat{\square}_h f ( a_k ) |\\
\leq \frac{1}{ \| w (  a_k  - \cdot ) \|_{L^1 ( B_h ( a_k ) \setminus B_\delta (a_k) ) }} \left| \int_{\sigma_k} \frac{ f ( a_k ) - f (y) }{ |  a_k  - y|^2 } w ( a_k -y) { \ }d y \right|\\
+ \frac{ \| w (  a_k  - \cdot )\|_{L^1 ( \sigma_k )}  }{ \| w ( a_k  - \cdot ) \|_{L^1 ( B_h ( a_k ) \setminus \sigma_k )} \| w \|_{L^1 (B_{h, \delta }) }} \left| \sum_{i \in \mathcal{R} (a_k , h)} \int_{\sigma_i} \frac{ f ( a_k ) - f (y) }{ |  a_k  - y|^2 } w ( a_k -y) { \ }d y \right|\\
\leq 2 \frac{ \| w (  a_k  - \cdot )/ | a_k - \cdot | \|_{L^1 ( \sigma_k \setminus B_\delta ( a_k ) )}}{\| w ( a_k  - \cdot ) \|_{L^1 ( B_h ( a_k ) \setminus B_\delta (a_k) )}} | f |_{C^1}.
\end{multline*}
Therefore, we have \eqref{eq62}.

Next, we derive \eqref{eq63}. A direct calculation gives
\begin{equation}\label{eq69}
\frac{1}{4} \{ \widehat{\square}_h f ( a_k ) - \breve{\square}_h f ( a_k ) \} = O_4 ( a_k ) + O_5 ( a_k ) + O_6 ( a_k ) + O_7 ( a_k ).
\end{equation}
Here
\begin{equation*}
O_4 ( a_k ) : = \frac{ \displaystyle{ \sum_{i \in \mathcal{R} ( a_k ,h) }\int_{\sigma_i} \frac{ f ( a_k ) - f (y)}{ |  a_k  - y | } \left( \frac{1}{ |  a_k  - y | } - \frac{1}{ |  a_k  - a_i | } \right) w ( a_k -y) { \ }d y }  }{ \displaystyle{ \sum_{j \in \mathcal{R} ( a_k ,h) }\int_{\sigma_j} w ( a_k - z ) { \ }d z }  },
\end{equation*}
\begin{equation*}
O_5 ( a_k ) := \frac{ \displaystyle{ \sum_{i \in \mathcal{R} ( a_k ,h) }\int_{\sigma_i} \frac{ f (a_i) - f (y) }{ |  a_k  - y | |  a_k  - a_i| } w ( a_k -y) { \ }d y }  }{ \displaystyle{ \sum_{j \in \mathcal{R} ( a_k ,h) }\int_{\sigma_j} w ( a_k - z ) { \ }d z }  } ,
\end{equation*}
\begin{equation*}
O_6 ( a_k ) := \frac{ \displaystyle{ \sum_{i \in \mathcal{R} ( a_k ,h) }\int_{\sigma_i} \frac{ f ( a_k ) - f (a_i) }{ |  a_k  - y | } \left( \frac{1}{ |  a_k  - y | } - \frac{1}{ |  a_k  - a_i | } \right) w ( a_k -y) { \ }d y }  }{ \displaystyle{ \sum_{j \in \mathcal{R} ( a_k ,h) }\int_{\sigma_j} w ( a_k - z ) { \ }d z }  } ,
\end{equation*}
and
\begin{equation*}
O_7 ( a_k ) := \frac{ \displaystyle{ \sum_{i \in \mathcal{R} ( a_k ,h) }\int_{\sigma_i} \frac{ f ( a_k ) - f (a_i) }{ |  a_k  - a_i | } \left( \frac{1}{ |  a_k  - y | } - \frac{1}{ |  a_k  - a_i | } \right) w ( a_k -y) { \ }d y }  }{ \displaystyle{ \sum_{j \in \mathcal{R} ( a_k ,h) }\int_{\sigma_j} w ( a_k - z ) { \ }d z }  } .
\end{equation*}
Since
\begin{equation*}
\left| \frac{1}{ |  a_k  - y|} - \frac{1}{ |  a_k - a_i |} \right| \leq \frac{ | a_i - y |}{ |  a_k  - y| |  a_k  - a_i| }, 
\end{equation*}
we use the mean-value theorem and \eqref{Eq213} to check that
\begin{align}
\sum_{p=4}^7 | O_p ( a_k ) | & \leq 4| f |_{C^1}  \frac{ \displaystyle{ \sum_{i \in \mathcal{R} ( a_k ,h) }\int_{\sigma_i} \frac{ | y - a_i| }{ |  a_k  - a_i |} \frac{w ( a_k -y)}{ |  a_k  - y | } { \ }d y }  }{ \displaystyle{ \sum_{j \in \mathcal{R} ( a_k ,h) }\int_{\sigma_j} w ( a_k -y) { \ }d y }  }\notag\\
& \leq 4 \bigg( \frac{r_\sigma}{\lambda h} \frac{ \left\| \frac{w (  a_k  - \cdot )}{ | a_k - \cdot | }\right\|_{L^1 ( B_h (a_k) \setminus \sigma_k ) }  }{ \displaystyle{ \| w (  a_k  - \cdot ) \|_{L^1 (B_h ( a_k ) \setminus \sigma_k )} }  }
 + \frac{ \left\| \frac{w (  a_k  - \cdot )}{| a_k - \cdot |} \right\|_{L^1 ( B_{\lambda h + r_\sigma} (a_k) \setminus \sigma_k ) }  }{ \displaystyle{ \| w (  a_k  - \cdot ) \|_{L^1 (B_h ( a_k ) \setminus \sigma_k )} }  } \bigg) | f |_{C^1} .\label{Eq610}
\end{align}
From \eqref{eq69} and \eqref{Eq610}, we have \eqref{eq63}.

Finally, we derive \eqref{eq64}. Using \eqref{eq25} and the mean-value theorem to see that
\begin{multline*}
\frac{1}{4} | \breve{\square}_h f ( a_k ) - \widetilde{\square}_h f ( a_k ) |\\
 \leq \frac{ \pi L_w }{ \| w (  a_k  - \cdot ) \|_{L^1 ( B_h ( a_k ) \setminus \sigma_k )} }
 \left( \frac{ 2 r_\sigma h}{ \lambda} + ( \lambda h + r_\sigma )^2 + \frac{ \displaystyle{\sum_{i \in \mathcal{R} (a_k , \lambda h )} V_i \frac{w (a_k - a_i)}{|a_k - a_i|} } }{ \displaystyle{\sum_{j \in \mathcal{R} (a_k , h )} V_j w (a_k - a_j)  } }  \right) | f |_{C^1}.
\end{multline*}
This is \eqref{eq64}. Therefore, the lemma follows.
\end{proof}

Finally, we prove Theorem \ref{thm16}.
\begin{proof}(Theorem \ref{thm16})
Using Lemmas \ref{lem61} and \ref{lem62}, we prove Theorem \ref{thm16}.
\end{proof}

\section{Appendix: Applications of Main Results}\label{sect7}
We state an application of the main results of this paper. We consider the following case:
\begin{equation*}
w (x) = 
\begin{cases}
1, & x \in \overline{B}_{h, \delta},\\
0, & x \in \mathbb{R}^2 \setminus \overline{B}_{h, \delta}.
\end{cases}
\end{equation*}
It is easy to check that $L_w = 0$ and that for each $0 <  q < p$,
\begin{align*}
\| w ( a_k - \cdot ) \|_{L^1 ( B_p (a_k) \setminus B_q ( a_k )  )} & = \pi ( p^2 - q^2 ),\\
\| w ( a_k - \cdot )/|a_k - \cdot | \|_{L^1 ( B_p (a_k) \setminus B_q ( a_k )  )} & = 2 \pi ( p - q ),\\
\| | a_k - \cdot | w ( a_k - \cdot ) \|_{L^1 ( B_p (a_k) \setminus B_q ( a_k )  )} & = \frac{ 2 \pi }{3} ( p^3 - q^3 ),\\
\| | a_k - \cdot |^2 w ( a_k - \cdot ) \|_{L^1 ( B_p (a_k) \setminus B_q ( a_k )  )} & = \frac{ \pi }{2} ( p^4 - q^4 ).
\end{align*}
In this section, we assume that $\delta = r_\sigma / 2$, 
\begin{align*}
\| w ( x - \cdot ) \|_{L^1 ( B_h (a_k) \setminus \sigma_k )} & = \| w ( x - \cdot ) \|_{L^1 ( B_{h,\delta} (a_k) )},\\
\frac{ \displaystyle{\sum_{i \in \mathcal{R} (a_k , h )} V_i | a_k - a_i | w (a_k - a_i)}}{ \displaystyle{\sum_{j \in \mathcal{R} (a_k , h )} V_j | a_k - a_j |^2 w (a_k - a_j)  } } & = \frac{ \| | a_k - \cdot | w (a_k - \cdot ) \|_{L^1 (B_{h} (a_k) \setminus B_\delta (a_k) )}}{ \| | a_k - \cdot |^2 w (a_k - \cdot ) \|_{L^1 (B_{h} (a_k) \setminus B_\delta (a_k) )}} .
\end{align*}
Suppose that there is $C_* > 1$ such that
\begin{equation*}
h = C_* r_\sigma .
\end{equation*}
Assume that $\mathcal{R}( a_k, \lambda h ) \neq \emptyset$ for some $0 < \lambda < 1$. Using Theorems \ref{thm13}-\ref{thm16} and $\| \cdot \|_{L^1 (\sigma_k \setminus B_\delta (a_k))} \leq \| \cdot \|_{L^1 (B_{r_\sigma} (a_k) \setminus B_\delta (a_k) )}$, we have the following corollary.
\begin{corollaryz}\label{cor71}{ \ }For each $f \in C^3 (\overline{\Omega}_H)$,
\begin{align}
| f (a_k) - \widetilde{\Pi}_h f (a_k) | & \leq ( C_* r_\sigma + r_\sigma )| f |_{C^1} + \frac{3/2}{C_*^2 - 1/4} | f |_{C^0},\\
| \nabla f (a_k) - \widetilde{\nabla}_h f (a_k) | & \leq 4 C_* r_\sigma | f |_{C^2} + \bigg( \frac{ 8}{  \lambda C_*} + \frac{8 ( \lambda C_* + 1)^2 + 1}{ C_*^2 - 1/4  } \bigg) | f |_{C^1},\label{eq72}\\
| \Delta f (a_k) - \widetilde{\square}_h f (a_k) | & \leq 24 C_* r_\sigma | f |_{C^3} + \bigg( \frac{32 \lambda C_* + 24}{r_\sigma ( C_*^2 - 1/4)} + \frac{32}{\lambda C_* r_\sigma ( C_* + 1/2)}  \bigg) | f |_{C^1},
\end{align}
\begin{multline}
| \Delta f (a_k) - \widetilde{\Delta}_h f (a_k) | \leq 24 C_* r_\sigma | f |_{C^3} + \frac{1}{3 r_\sigma ( C_*^4 - 1/{16})} \cdot\\
\bigg( 8 + 24 (\lambda C_* + 1 )^2 + \frac{ 5 6 C_*^3 - 7}{ 3( C_*^4 - 1/{16})} + \frac{16 C_*^3 - 2}{\lambda C_*} + \frac{ C_* ( 6 4 C_*^3 - 8 )}{C_*^2+ 1/4}  \bigg) | f |_{C^1}.
\end{multline}
Moreover, the following two assertions hold:\\
$(\mathrm{i})$ If $r_\sigma = 10^{-15}$, $C_* = 10^{12}$, $\lambda = 10^{-6}$, then
\begin{align}
| f (a_k) - \widetilde{\Pi}_h f (a_k) | & \leq \frac{1}{999} | f |_{C^1} + \frac{1}{10^{23}} | f |_{C^0} ,\\
| \nabla f (a_k) - \widetilde{\nabla}_h f (a_k) | & \leq \frac{1}{250} | f |_{C^2} + \frac{1}{10^5} | f |_{C^1},\label{eq76}\\
| \Delta f (a_k) - \widetilde{\Delta}_h f (a_k) | & \leq \frac{3}{125} | f |_{C^3} + \frac{1}{100} | f |_{C^1},\\
| \Delta f (a_k) - \widetilde{\square}_h f (a_k) | & \leq \frac{3}{125} | f |_{C^3} + \frac{1}{10} | f |_{C^1}.
\end{align}
$(\mathrm{ii})$ If $r_\sigma = 10^{-2}$, $C_* = 4$, $\lambda = 1/2$, then
\begin{align}
| f (a_k) - \widetilde{\Pi}_h f (a_k) | & \leq \frac{1}{20} | f |_{C^1} + \frac{1}{10} |f|_{C^0},\\
| \nabla f (a_k) - \widetilde{\nabla}_h f (a_k) | & \leq \frac{4}{25} | f |_{C^2} + 10 | f |_{C^1},\\
| \Delta f (a_k) - \widetilde{\Delta}_h f (a_k) | & \leq \frac{24}{25} | f |_{C^3} + 300 | f |_{C^1},\\
| \Delta f (a_k) - \widetilde{\square}_h f (a_k) | & \leq \frac{24}{25} | f |_{C^3} + 1000 | f |_{C^1}.
\end{align}
\end{corollaryz}

\begin{proof}(Corollary \ref{cor71}) We only show \eqref{eq72} and \eqref{eq76} since other cases are similar. Fix $f \in ^3 (\overline{\Omega}_H)$.

We first show \eqref{eq72}. Since $L_w =0$, it follows from Theorem \ref{thm14} to see that
\begin{multline}\label{eq713}
| \nabla f ( a_k ) - \widetilde{\nabla}_h f (a_k ) | \leq 4 h | f |_{C^2}\\ + \left\{ \frac{8 r_\sigma}{ \lambda h} + 4 \frac{ \| w (a_k - \cdot ) \|_{L^1 (\sigma_k \setminus B_\delta (a_k) )}}{ \| w ( a_k - \cdot ) \|_{L^1 (B_h (a_k) \setminus B_\delta (a_k) )}} + 8 \frac{ \displaystyle{ \| w ( a_k - \cdot ) \|_{L^1 ( B_{ \lambda h + r_\sigma } ( a_k ) \setminus \sigma_k )} } }{\| w ( a_k - \cdot ) \|_{L^1 ( B_h ( a_k ) \setminus \sigma_k )}} \right\} | f |_{C^1}.
\end{multline}
By the assumptions of Section \ref{sect7}, we check that
\begin{multline*}
\text{(R.H.S.) of } \eqref{eq713}\\
\leq 4 C_* r_\sigma | f |_{C^2} + \left( \frac{8 r_\sigma}{ \lambda C_* r_\sigma} + \frac{4 \pi (r_\sigma^2 - r_\sigma^2/4)}{\pi (C_*^2 r_\sigma^2 - r_\sigma^2/4)} + \frac{8 \pi \{ (\lambda C_* r_\sigma + r_\sigma)^2 - r_\sigma^2/4 \} }{\pi (C_*^2 r_\sigma^2 - r_\sigma^2/4)} \right) | f |_{C^1}\\
 = 4 C_* r_\sigma | f |_{C^2} + \bigg( \frac{ 8}{ \lambda C_*} + \frac{8 ( \lambda C_* + 1)^2 + 1}{ C_*^2 - 1/4 } \bigg) | f |_{C^1}.
\end{multline*}
Therefore, we obtain \eqref{eq72}.

Next we show \eqref{eq76}. A direct calculation shows that
\begin{align*}
4 C_* r_\sigma = 4 \cdot 10^{12} \cdot 10^{-15} = \frac{4}{1000} = \frac{1}{250},
\end{align*}
and that
\begin{align*}
\frac{8}{\lambda C_*} + \frac{8(\lambda C_* + 1)^2+ 1}{C_*^2 - 1/4} & = \frac{8}{10^6} + \frac{8 (10^6 + 1)^2 + 1 }{10^{24} -1}\\
& \leq \frac{8}{10^6} + \frac{1}{10^6} \leq \frac{1}{10^5}.
\end{align*}
From \eqref{eq72}, we have \eqref{eq76}. Therefore, Corollary \ref{cor71} is proved.
\end{proof}

\end{document}